\documentclass[a4paper, 11pt]{article}
\usepackage[T2A,T1]{fontenc}
\usepackage[utf8]{inputenc}
\usepackage[german, russian, english]{babel}
\usepackage{amsmath, amsthm, amssymb,amsfonts}
\usepackage{mathtools}
\usepackage{mathrsfs}
\usepackage{enumerate}
\usepackage{graphicx}
\usepackage{color}
\usepackage{hyperref}

\numberwithin{equation}{section}

\theoremstyle{definition}
\newtheorem{theorem}{Theorem}
\newtheorem{corollary}[theorem]{Corollary}

\newtheorem{remark}[theorem]{Remark}
\newtheorem{definition}[theorem]{Definition}

\numberwithin{theorem}{section}

\newcommand{\NN}{\mathbb{N}}

\begin{document}

\title{Representation of Classical Solutions to a Linear Wave Equation with Pure Delay}

\author{Denys Ya. Khusainov\footnote{Department of Cybernetics, Kyiv National Taras Shevchenko University, Ukraine},
Michael Pokojovy\footnote{Department of Mathematics and Statistics, University of Konstanz, Germany},
Elvin I. Azizbayov\footnote{Department of Mechanics and Mathematics, Baku State University, Azerbaijan}}

\date{\today}

\maketitle

\begin{abstract}
    For a wave equation with pure delay, we study an inhomogeneous initial-boundary value problem in a bounded 1D domain.
    Under smoothness assumptions, we prove unique existence of classical solutions for any given finite time horizon
    and give their explicit representation.
    Continuous dependence on the data in a weak extrapolated norm is also shown.
\end{abstract}

\section{Introduction}
    The wave equation is a typical linear hyperbolic second-order partial differential equation
    which naturally arises when modeling phenomena of continuum mechanics such as
    sound, light, water or other kind of waves in acoustics, (electro)magnetics, elasticity and fluid dynamics, etc.
    (cf. \cite{EckGaKna2009, Ti1990}).
    Providing a rather adequate description of physical processes,
    partial differential equations, or equations with distributed parameters in general,
    have found numerous applications in mechanics, medicine, ecology, etc.
    Introducing after-effects such as delay into these equations has gained a lot of attention
    over several past decades. See e.g., \cite{AzKhu2012, BaPia2005, ElNo1973, Ha1977}.
    Mathematical treatment of such systems requires additional carefulness since distributed systems
    with delay often turn out to be even ill-posed (cf. \cite{Da1997, DreQuiRa2009, Ra2012}).

    In the present paper, we consider an initial-boundary value problem for a general linear wave equation
    with pure delay and constant coefficients in a bounded interval subject to non-homogeneous Dirichlet boundary conditions.
    To solve the equation, we employ Fourier's separation method as well as the special functions
    referred to as delay sine and cosine functions
    which were introduced in \cite{KhuDiRuLu2008,KhuIvKo2006}.
    We prove the existence of a unique classical solution on any finite time interval,
    show its continuous dependence on the data in a very weak extrapolated norm,
    give its representation as a Fourier series and prove its absolute and uniform convergence
    under certain conditions on the data.

\section{Equation with pure delay}
    For $T > 0$, $l > 0$, we consider the following linear wave equation in a bounded interval $(0, l)$ with a single delay
    being a second order partial difference-differential equation for an unknown function $\eta$
    \begin{equation}
        \partial_{tt} \eta(t, x) = a^{2} \partial_{xx} \eta(t - \tau, x) + b \partial_{x} \eta(t - \tau, x) + d \eta(t - \tau, x) + g(t, x)
        \text{ for } (t, x) \in (0, T) \times (0, l)
        \label{EQ1_1}
    \end{equation}
    subject to non-homogeneous Dirichlet boundary conditions and initial conditions
    \begin{equation}
        \begin{split}
            \eta(t, 0) = \theta_{1}(t), \quad \eta(t, l) = \theta_{2}(t) \text{ for } t > -\tau, \\
            \eta(t, x) = \psi(t, x) \text{ for } (t, x) \in (-\tau, 0) \times (0, l).
        \end{split}
        \label{EQ1_2}
    \end{equation}
    Since we are interested in studying classical solutions,
    the following compatibility conditions are required to assure for smoothness of the solution on the boundary of time-space cylinder
    \begin{equation}
        \psi(t, 0) = \theta_{1}(t), \quad \psi(t, l) = \theta_{2}(t) \text{ for } t > -\tau. \notag
    \end{equation}

    \begin{definition}
        Under a classical solution to the problem (\ref{EQ1_1}), (\ref{EQ1_2}) we understand
        a function $\eta \in \mathcal{C}^{0}\big([-\tau, T] \times [0, l]\big)$ which satisfies
        $\partial_{tt} \eta, \partial_{tx} \eta, \partial_{xx} \eta \in \mathcal{C}^{0}\big([-\tau, 0] \times [0, l]\big)$
        as well as
        $\partial_{tt} \eta, \partial_{tx} \eta, \partial_{xx} \eta \in \mathcal{C}^{0}\big([0, T] \times [0, l]\big)$
        and, being plugged into Equations (\ref{EQ1_1}), (\ref{EQ1_2}), turns them into identity.
    \end{definition}

    \begin{remark}
        The previous definition does not impose any continuity of time derivatives in $t = 0$.
        If the continuity or even smoothness are desired,
        additional compatibility conditions on the data, including $g$, are required.
    \end{remark}

    Let $\|\cdot\|_{k, 2} := \|\cdot\|_{H^{k, 2}((0, l))}$, $k \in \NN_{0}$, denote the standard Sobolev norm (cf. \cite{Ad1975})
    and $\|\cdot\|_{-k, 2} := \|\cdot\|_{H^{-k, 2}((0, l))}$ denote the norm of corresponding negative Sobolev space.
    We introduce the norm $\|\cdot\|_{X} := \sqrt{\sum\limits_{k = 0}^{\infty} \|\cdot\|_{-k, 2}^{2}}$
    and define the Hilbert space $X$ as a completion of $L^{2}\big((0, l)\big)$ with respect to $\|\cdot\|_{X}$.
    Obviously, $X \hookrightarrow \left(\mathcal{D}\big((0, l)\big)\right)'$,
    i.e., $X$ can be continuously embedded into the space of distributions.

    With this notation, we easily see that $\mathcal{A} := a^{2} \partial_{x}^{2} + b \partial_{x} + d$
    (with $\partial_{x}$ denoting the distributional derivative)
    is a bounded linear operator on $X$ since
    \begin{equation}
        \begin{split}
            \|\mathcal{A}\|_{L(X)} &= \sup_{\|u\|_{X} = 1} \|\mathcal{A} u\|_{X} =
            \sup_{\|u\|_{X} = 1} \sqrt{\sum_{k = 0}^{\infty} \|a^{2} \partial_{x}^{2} u + b \partial_{x} u + du\|_{-k, 2}^{2}} \\
            &\leq \sup_{\|u\|_{X} = 1} \sum_{k = 0}^{\infty} \left(a^{2} \|u\|_{-k-2, 2} + b \|u\|_{-k-1, 2} + d \|u\|_{-k, 2}\right) \\
            &\leq \sup_{\|u\|_{X} = 1} (a^{2} + b + c) \|u\|_{X} = a^{2} + b + d.
        \end{split}
        \notag
    \end{equation}

    \begin{theorem}
        \label{THEOREM1_1}
        There exists a constant $C > 0$, dependent only on $a, b, d, l, \tau, T$,
        such that the estimate
        \begin{equation}
            \begin{split}
                \max_{t \in [0, T]} \left(\|\eta(\cdot, t)\|_{X}^{2} + \|\eta_{t}(\cdot, t)\|_{X}^{2}\right) &\leq
                C \left(\|\psi(0, \cdot)\|_{X}^{2} + \|\psi_{t}(0, \cdot)\|_{X}^{2} \right) +
                C \int_{-\tau}^{0} \left(\|\psi(s, \cdot)\|_{X}^{2} + \|\psi_{t}(s, \cdot)\|_{X}^{2}\right) \mathrm{d}t + \\
                &\phantom{=\;\;} C \int_{0}^{T} \left(\|g(s, \cdot)\|_{X}^{2} + |\theta_{1}(s)|^{2} + |\theta_{2}(s)|^{2}\right) \mathrm{d}s
            \end{split} \notag
        \end{equation}
    \end{theorem}
    holds true for any classical solution of Equations (\ref{EQ1_1}), (\ref{EQ1_2}).
    \begin{proof}
        Let $\eta$ be the classical solution to Equations (\ref{EQ1_1}), (\ref{EQ1_2}). We define
        \begin{equation}
            w(t, x) :=
            \begin{cases}
                \eta(t, x) & \text{ for } (t, x) \in [-\tau, 0] \times [0, l], \\
                \eta(t, x) - \theta_{1}(t) - \tfrac{x}{l}\left(\theta_{2}(t) - \theta_{1}(t)\right) &
                \text{ for } (t, x) \in (0, T] \times [0, l]. \\
            \end{cases}
            \notag
        \end{equation}
        Then $w$ satisfies homogeneous Dirichlet boundary conditions and solves the equation
        \begin{equation}
            \partial_{tt} w(t, \cdot) = \mathcal{A} w(t - \tau, \cdot) + f(t, \cdot)
            \label{EQ1_3}
        \end{equation}
        in the extrapolated space $X$ with
        \begin{equation}
            f(t, \cdot) = g(t, \cdot) + b\left(\theta_{2}(t) - \theta_{1}(t)\right) + \theta_{1}(t)
            + \tfrac{dx}{l} \left(\theta_{2}(t) - \theta_{1}(t)\right).
            \notag
        \end{equation}
        We multiply the equation with $w_{t}(t, \cdot)$ in the scalar product of $X$ and use Young's inequality to get the estimate
        \begin{equation}
            \begin{split}
                \partial_{t} \|w_{t}(t, \cdot)\|_{X}^{2} &=
                \langle Aw(t - \tau, \cdot), w_{t}(t, \cdot)\rangle_{X} + \langle f(t, \cdot), w_{t}(t, \cdot)\rangle_{X} \\
                &\leq \|w(t - \tau, \cdot)\|_{X}^{2} + \left(1 + \|\mathcal{A}\|_{L(X)}^{2}\right) \|w_{t}(t, \cdot)\|_{X}^{2} +
                \|f(t, \cdot)\|_{X}^{2}.
            \end{split}
            \label{EQ1_4}
        \end{equation}
        As in \cite{KhuPoAz2013}, we introduce the history variable
        \begin{equation}
            z(s, t, x) := w(t - s, x) \text{ for } (s, t, x) \in [0, \tau] \times [0, T] \times [0, l] \notag
        \end{equation}
        and obtain
        \begin{equation}
            z_{t}(s, t, x) + z_{s}(s, t, x) = 0 \text{ for } (s, t, x) \in (0, \tau) \times (0, T) \times (0, l). \notag
        \end{equation}
        Multiplying these identities with $w(t, \cdot)$ in $X$ and performing a partial integration, we find
        \begin{equation}
            \partial_{t} \int_{0}^{\tau} \|z(s, t, \cdot)\|_{X}^{2} \mathrm{d}s =
            -\int_{0}^{\tau} \partial_{s} \|z(s, t, \cdot)\|_{X}^{2} \mathrm{d}s =
            \|w(t, \cdot)\|_{X}^{2} - \|w(t - \tau, \cdot)\|_{X}^{2}.
            \label{EQ1_5}
        \end{equation}
        Adding equations (\ref{EQ1_4}) and (\ref{EQ1_5}) to the trivial identity
        \begin{equation}
            \partial_{t} \|w(t, \cdot)\|_{X}^{2} \leq \|w(t, \cdot)\|_{X}^{2} + \|w_{t}(t, \cdot)\|_{X}^{2}, \notag
        \end{equation}
        we obtain
        \begin{equation}
            \partial_{t} \left(\|w(t, \cdot)\|_{X}^{2} + \|w_{t}(t, \cdot)\|_{X}^{2} + \int_{-\tau}^{0} \|w(s, t, \cdot)\|_{X}^{2} \mathrm{d}s\right) \leq
            \left(2 + \|\mathcal{A}\|_{L(X)}\right) \|w(t, \cdot)\|_{X}^{2} + \|w_{t}(t, \cdot)\|_{X}^{2} + \|f(t, \cdot)\|_{X}^{2}. \notag
        \end{equation}
        Thus, we have shown
        \begin{equation}
            \partial_{t} E(t) \leq \left(2 + \|\mathcal{A}\|_{L(X)}\right) E(t) + \|f(t, \cdot)\|_{X}^{2},
            \label{EQ1_6}
        \end{equation}
        where
        \begin{equation}
            E(t) := \|w(t, \cdot)\|_{X}^{2} + \|w_{t}(t, \cdot)\|_{X}^{2} + \int_{-\tau}^{0} \|z(s, t, \cdot)\|_{X} \mathrm{d}s. \notag
        \end{equation}
        From Equation (\ref{EQ1_6}) we conclude
        \begin{equation}
            E(t) \leq E(0) + \left(2 + \|\mathcal{A}\|_{L(X)}\right) \int_{0}^{t} E(s) \mathrm{d}s + \int_{0}^{t} \|f(t, \cdot)\|_{X}^{2} \mathrm{d}s.
            \notag
        \end{equation}
        Using now the integral form of Gronwall's inequality, we obtain
        \begin{equation}
            \begin{split}
                E(t) &\leq E(0) + \int_{0}^{t} \|f(t, \cdot)\|_{X}^{2} \mathrm{d}s +
                \int_{0}^{t} e^{\left(2 + \|\mathcal{A}\|_{L(X)}\right) (t - s)} \left(E(0) + \int_{0}^{s} \|f(\xi, \cdot)\|_{X}^{2} \mathrm{d}\xi\right) \mathrm{d}s \\
                &\leq \tilde{C} \left(E(0) + \int_{0}^{T} \|f(s, \cdot)\|_{X}^{2} \mathrm{d}s\right)
            \end{split}
            \label{EQ1_7}
        \end{equation}
        for certain $\tilde{C} > 0$. Taking into account
        \begin{equation}
            \begin{split}
                c_{1} \left(\|w(t, \cdot)\|_{X}^{2} + |\theta_{1}(t)|^{2} + |\theta_{2}(t)|^{2}\right) &\leq \|\eta(t, \cdot)\|_{X}^{2} \leq
                C_{1} \left(\|w(t, \cdot)\|_{X}^{2} + |\theta_{1}(t)|^{2} + |\theta_{2}(t)|^{2}\right), \\
                c_{2} \left(\|f(t, \cdot)\|_{X}^{2} + |\theta_{1}(t)|^{2} + |\theta_{2}(t)|^{2}\right) &\leq \|g(t, \cdot)\|_{X}^{2} \leq
                C_{2} \left(\|f(t, \cdot)\|_{X}^{2} + |\theta_{1}(t)|^{2} + |\theta_{2}(t)|^{2}\right)
            \end{split}
            \notag
        \end{equation}
        for some constants $c_{1}, c_{2}, C_{1}, C_{2} > 0$ and exploiting the definition of $E(t)$,
        the proof is a direct consequence of Equation (\ref{EQ1_7}).
    \end{proof}

    \begin{corollary}
        Solutions of Equation (\ref{EQ1_1}), (\ref{EQ1_2}) are unique.
        The solution map
        \begin{equation}
            (\psi, g, \theta_{1}, \theta_{2}) \mapsto \eta \notag
        \end{equation}
        is well-defined, linear and continuous in the norms from Theorem \ref{THEOREM1_1}.
    \end{corollary}

    \begin{remark}
        It was essential to consider the weak space $X$.
        If the space corresponding to the usual wave equation is used, i.e., $(\eta, \eta_{t}) \in H^{1}_{0}\big((0, l)\big) \times L^{2}\big((0, l)\big)$,
        there follows from \cite{DreQuiRa2009}
        that Equation (\ref{EQ1_1}), (\ref{EQ1_2}) is an ill-posed problem due to the lack of continuous dependence
        on the data even in the homogeneous case.
    \end{remark}

    Next, we want to establish conditions on the data allowing for the existence of a classical solution.
    Performing the substitution
    \begin{equation}
        \xi(t, x) := e^{-\tfrac{b}{2a^{2}} x} \eta(t, x) \text{ for } (t, x) \in [-\tau, T] \times [0, l]
        \label{EQ1_8}
    \end{equation}
    with a new unknown function $\xi$ (cp. \cite{KhuPoAz2013}), the initial boundary value problem (\ref{EQ1_1}), (\ref{EQ1_2})
    can be written in the following simplified form with a self-adjoint operator on the right-hand side
    \begin{equation}
        \partial_{tt} \xi(t, x) = a^{2} \partial_{xx} \xi(t - \tau, x) + c \xi(t - \tau, x) + f(t, x) \text{ for } (t, x) \in (0, T) \times (0, l)
        \label{EQ1_9}
    \end{equation}
    with $c := d - \frac{b^{2}}{4 a^{2}}$ complemented by the following boundary and initial conditions
    \begin{align}
        \xi(t, 0) &= \mu_{1}(t), \xi(t, l) = \mu_{2}(t) \text{ for } t > -\tau
        \text{ with } \mu_{1}(t) := \theta_{1}(t), \quad \mu_{2}(t) := e^{\frac{b}{2a^{2}} l} \theta_{2},
        \label{EQ1_10} \\
        \xi(t, x) &= \varphi(t, x) \text{ for } (t, x) \in (-\tau, 0) \times (0, l) \text{ with } \varphi(t, x) := e^{\frac{b}{2a^{2} x}} \psi(t, x)
        \label{EQ1_11}
    \end{align}
    and
    \begin{equation}
        f(t, x) := e^{\frac{b}{2a^{2}} x} g(t, x) \text{ for } (t, x) \in [0, T] \times [0, l]. \notag
    \end{equation}

    The solution will be determined in the form
    \begin{equation}
        \xi(t, x) = \xi_{0}(t, x) + \xi_{1}(t, x) + G(t, x). \notag
    \end{equation}
    Here, $G$ is an arbitrary function with $\partial_{tt} G, \partial_{tx} G, \partial_{xx} G \in \mathcal{C}^{0}\big([-\tau, T] \times [0, l]\big)$
    satisfying the boundary conditions
    \begin{equation}
        G(t, 0) = \mu_{1}(t), \quad G(t, l) = \mu_{2}(t). \notag
    \end{equation}
    Assuming $\mu_{1}, \mu_{2} \in \mathcal{C}^{2}\big([-\tau, T]\big)$, we let
    \begin{equation}
        G(t, x) := \mu_{1}(t) + \tfrac{x}{l} \left(\mu_{2}(t) - \mu_{1}(t)\right) \text{ for } (t, x) \in [-\tau, T] \times [0, l].
        \label{EQ1_12}
    \end{equation}
    \begin{itemize}
        \item $\xi_{0}$ solves the homogeneous equation
        \begin{equation}
            \partial_{tt} \xi_{0}(t, x) = a^{2} \partial_{xx} \xi_{0}(t - \tau, x) + c \xi_{0}(t - \tau, x)
            \label{EQ1_13}
        \end{equation}
        subject to homogeneous boundary and non-homogeneous initial conditions
        \begin{equation}
            \begin{split}
                \xi_{0}(t, 0) &\equiv 0, \xi_{0}(t, l) = 0 \text{ in } (-\tau, T), \\
                \xi_{0}(t, x) &= \Phi(t, x) \text{ for } (t, x) \in (-\tau, 0) \times (0, l) \text{ with } \Phi(t, x) := \varphi(t, x) - G(t, x).
            \end{split}
            \label{EQ1_14}
        \end{equation}
        In particular, with the function $G$ selected as in Equation (\ref{EQ1_12}), we obtain
        \begin{equation}
            \Phi(t, x) = \varphi(t, x) - \mu_{1}(t) - \tfrac{x}{l} \left(\mu_{2}(t) - \mu_{1}(t)\right) \text{ for }
            (t, x) \in [-\tau, 0] \times [0, l].
            \label{EQ1_15}
        \end{equation}

        \item $\xi_{1}$ solves the non-homogeneous equation
        \begin{equation}
            \partial_{tt} \xi_{1}(t, x) = a^{2} \partial_{xx} \xi_{1}(t - \tau, x) + c \xi_{1}(t - \tau, x) + F(t, x)
            \text{ for } (t, x) \in (0, T) \times (0, l)
            \label{EQ1_16}
        \end{equation}
        with
        \begin{equation}
            F(t, x) := a^{2} \partial_{xx} G(t - \tau, x) + c G(t - \tau, x) - \partial_{tt} G(t, x) \text{ for } (t, x) \in [0, T] \times [0, l]
            \label{EQ1_17}
        \end{equation}
        subject to homogeneous boundary and initial conditions.
        For $G$ from Equation (\ref{EQ1_12}), we have
        \begin{equation}
            F(t, x) = f(t, x) + c \left(\mu_{1}(t - \tau) + \tfrac{x}{l} \left(\mu_{2}(t - \tau) - \mu_{1}(t - \tau)\right)\right) -
            \left(\ddot{\mu}_{1}(t) + \tfrac{x}{l} \left(\ddot{\mu}_{2}(t) - \ddot{\mu}_{1}(t)\right)\right).
            \label{EQ1_18}
        \end{equation}
    \end{itemize}

    \section{Homogeneous equation}
    In this section, we obtain a formal solution to the initial-boundary value problem (\ref{EQ1_13})
    with initial and boundary conditions given in Equations (\ref{EQ1_10}), (\ref{EQ1_11}).
    We exploit Fourier's separation method to determine $\xi_{0}$ in the product form $\xi_{0}(t, x) = T(t) X(x)$.
    After plugging this ansatz into Equation (\ref{EQ1_13}), we find
    \begin{equation}
        X(x) \ddot{T}(t) = a^{2} X''(x) T(t - \tau) + c X(x) T(t - \tau). \notag
    \end{equation}
    Hence,
    \begin{equation}
        X(x) \left(\ddot{T}(t) - c T(t - \tau)\right) = a^{2} X''(x) T(t - \tau). \notag
    \end{equation}
    By formally separating the variables, we deduce
    \begin{equation}
        \frac{X''(x)}{X(x)} = \frac{\ddot{T}(t) - c T(t - \tau)}{a^{2} T(t - \tau)} = -\lambda^{2}. \notag
    \end{equation}
    Thus, the equation can be decoupled as follows
    \begin{equation}
        \ddot{T}(t) + \left(a^{2} \lambda^{2} - c\right) T(t - \tau) = 0, \quad
        X''(x) + \lambda^{2} X(x) = 0.
        \label{EQ2_1}
    \end{equation}
    These are linear second order ordinary (delay) differential equations with constant coefficients.

    Due to the zero boundary conditions for $\xi_{0}$, the boundary conditions for the second equation in (\ref{EQ2_1}) will also be homogeneous, i.e.,
    \begin{equation}
        X(0) = 0, \quad X(l) = 0. \notag
    \end{equation}
    Therefore, we obtain a Sturm \& Liouville problem admitting non-trivial solutions only for the eigennumbers
    \begin{equation}
        \lambda^{2} = \lambda_{n}^{2} = \left(\tfrac{\pi n}{l}\right)^{n} \text{ for } n \in \NN \notag
    \end{equation}
    and the corresponding eigenfunctions
    \begin{equation}
        X_{n}(x) = \sin\left(\tfrac{\pi n}{l} x\right) \text{ for } n \in \NN. \notag
    \end{equation}
    Assuming
    \begin{equation}
        \left(\tfrac{\pi a}{l}\right)^{2} - c > 0, \notag
    \end{equation}
    we denote
    \begin{equation}
        \omega_{n} = \sqrt{\left(\tfrac{\pi n}{a}\right)^{2} - c} \text{ for } n \in \NN \notag
    \end{equation}
    and consider the first equation in (\ref{EQ2_1}), i.e.,
    \begin{equation}
        \ddot{T}(t) + \omega_{n}^{2} T(t - \tau) = 0 \text{ for } n \in \NN.
        \label{EQ2_2}
    \end{equation}
    The initial conditions for each of the equations in (\ref{EQ2_2}) can be obtained by expanding the initial data
    into a Fourier series with respect to the eigenfunction basis of the second equation in (\ref{EQ2_1})
    \begin{equation}
        \begin{split}
            \Phi(t, \cdot) &= \sum_{n = 1}^{\infty} \Phi_{n}(t) \sin\left(\tfrac{\pi n}{l} x\right)
            \text{ with }
            \Phi_{n}(t) = \frac{2}{l} \int_{0}^{l} \left(\varphi(t, s) - G(t, s)\right) \sin\left(\tfrac{\pi n}{l} s\right) \mathrm{d}s, \\
            \partial_{t} \Phi(t, \cdot) &= \sum_{n = 1}^{\infty} \dot{\Phi}_{n}(t) \sin\left(\tfrac{\pi n}{l} x\right)
            \text{ with }
            \dot{\Phi}_{n}(t) = \frac{2}{l} \int_{0}^{l} \left(\partial_{t} \varphi(t, s) - \partial_{t} G(t, s)\right) \sin\left(\tfrac{\pi n}{l} s\right) \mathrm{d}s
        \end{split}
        \label{EQ2_3}
    \end{equation}
    for $t \in [-\tau, T]$.
    Let us further determine the solution of the Cauchy problem associated with each of the equations in (\ref{EQ2_2})
    subject to the initial conditions from Equation (\ref{EQ2_3}).

    First, we briefly present some useful results from the theory of second order delay differential equations
    with pure delay obtained in \cite{KhuDiRuLu2008}.
    The authors considered a linear homogeneous second order delay differential equation
    \begin{equation}
        \ddot{x}(t) + \omega^{2} x(t - \tau) = 0 \text{ for } t \in (0, \infty), \quad
        x(t) = \beta(t) \text{ for } t \in [-\tau, 0].
        \label{EQ2_4}
    \end{equation}
    They introduced two special functions referred to as delay cosine and sine functions.
    Exploiting these functions, a unique solution to the initial value problem (\ref{EQ2_4}) was obtained.

    \begin{definition}
        Delay cosine is the function given as
        \begin{equation}
            \cos_{\tau}(\omega, t) =
            \left\{
            \begin{array}{cc}
                0, & -\infty < t < -\tau, \\
                1, & -\tau \leq t < 0, \\
                1 - \omega^{2} \tfrac{t^{2}}{2!}, & 0 \leq t < \tau, \\
                \vdots & \vdots \\
                1 - \omega^{2} \tfrac{t^{2}}{2!} + \omega^{4} \tfrac{(t - \tau)^{4}}{4!} - \dots +
                (-1)^{k} \omega^{2k} \tfrac{(t - (k - 1)\tau)^{2k}}{(2k)!}, &
                (k - 1) \tau \leq t < k\tau
            \end{array}\right.
            \label{EQ2_5}
        \end{equation}
        with $2k$-order polynomials on each of the intervals $(k - 1) \tau \leq t < k\tau$
        continuously adjusted at the nodes $t = k\tau$, $k \in \NN_{0}$.
    \end{definition}

    \begin{figure}[h!]
        \centering
        \includegraphics[scale = 0.5]{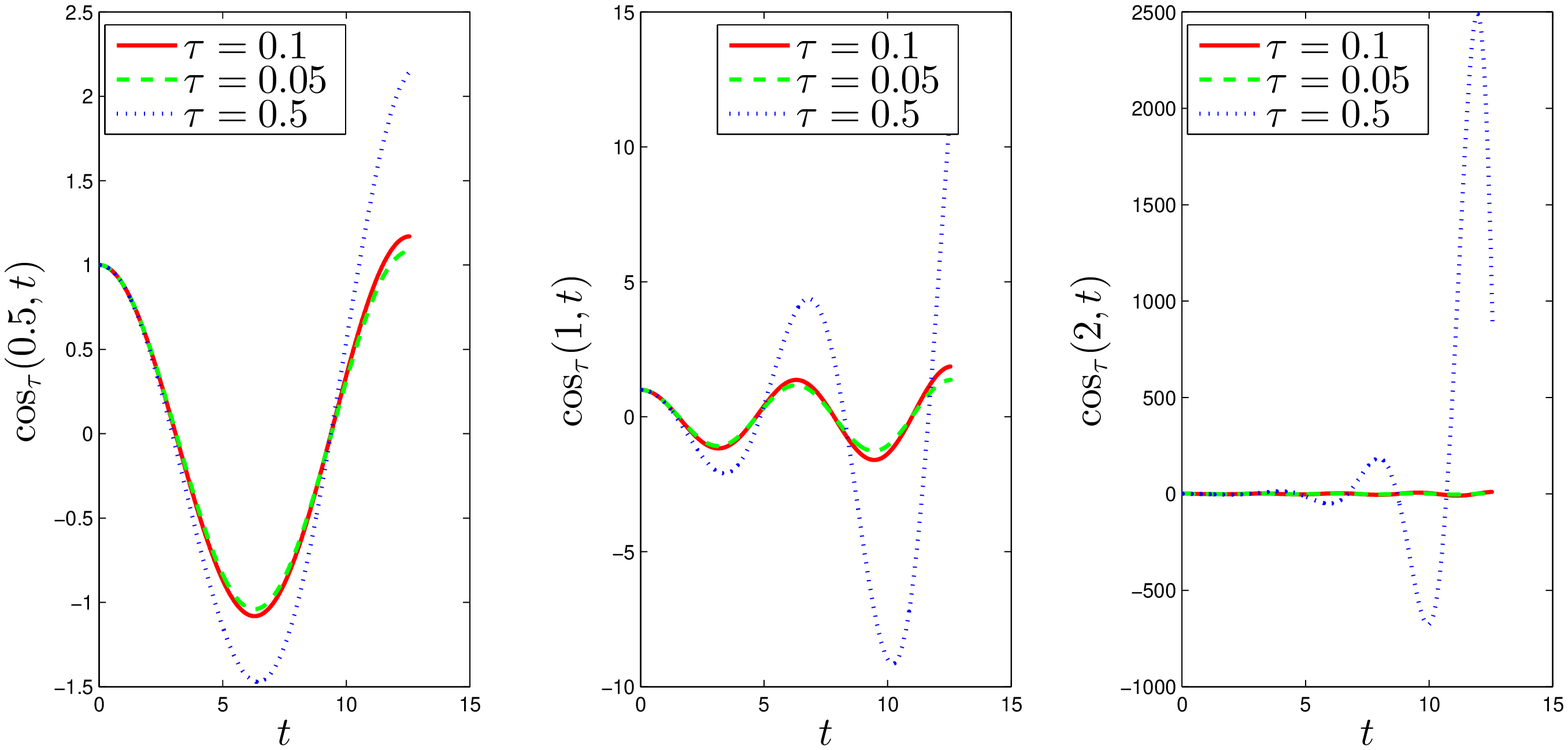}
        \caption{Delay cosine function}
    \end{figure}

    \begin{definition}
        Delay sine is the function given as
        \begin{equation}
            \cos_{\tau}(\omega, t) =
            \left\{
            \begin{array}{cc}
                0, & -\infty < t < -\tau, \\
                \omega (1 + \tau), & -\tau \leq t < 0, \\
                \omega (1 + \tau) - \omega^{3} \tfrac{t^{3}}{3!}, & 0 \leq t < \tau, \\
                \vdots & \vdots \\
                \omega (1 + \tau) - \omega^{3} \tfrac{t^{3}}{3!} + \dots +
                (-1)^{k} \omega^{2k+1} \tfrac{(t - (k - 1) \tau)^{2k + 1}}{(2k + 1)!}, &
                (k - 1) \tau \leq t < k\tau
            \end{array}\right.
            \label{EQ2_6}
        \end{equation}
        with $(2k+1)$-order polynomials on each of the intervals $(k - 1) \tau \leq t < k\tau$
        continuously adjusted at the nodes $t = k\tau$, $k \in \NN_{0}$.
    \end{definition}

    \begin{figure}[h!]
        \centering
        \includegraphics[scale = 0.5]{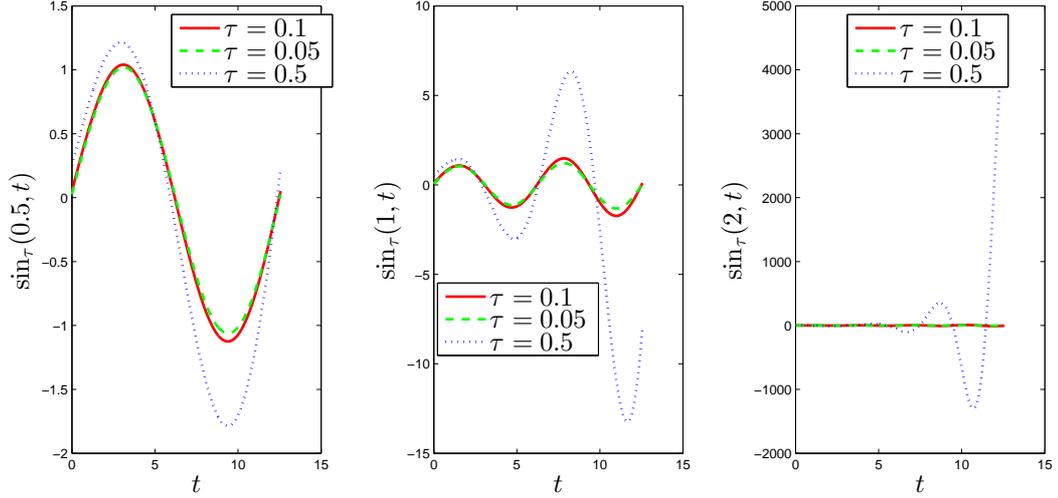}
        \caption{Delay sine function}
    \end{figure}

    There has further been proved that delay cosine uniquely solves the linear homogeneous second order ordinary delay differential equation
    with pure delay subject to the unit initial conditions $x \equiv 1$ in $[-\tau, 0]$
    and the delay sine in its turn solves Equation (\ref{EQ2_4}) subject to the initial conditions
    $x(t) = \omega(t + \tau)$ for $t \in [-\tau, 0]$.

    Using the fact above, the solution of the Cauchy problem was represented in the integral form.
    In particular, the solution $x$ to the homogeneous delay differential equation (\ref{EQ2_4})
    with the initial conditions $x \equiv \beta$ in $[-\tau, 0]$ for an arbitrary $\beta \in \mathcal{C}^{2}\left([-\tau, 0]\right)$
    was shown to be given as
    \begin{equation}
        x(t) = \beta(-\tau) \cos_{\tau}(\omega, t) + \tfrac{1}{\omega} \dot{\beta}(-\tau) \sin_{\tau}(\omega, t) +
        \tfrac{1}{\omega} \int_{-\tau}^{0} \sin_{\tau}(\omega, t - \tau - s) \ddot{\beta}(s) \mathrm{d}s.
        \label{EQ2_7}
    \end{equation}
    Turning back to the delay differential equation (\ref{EQ2_2}) with the initial conditions (\ref{EQ1_4}),
    we obtain their unique solution in the form
    \begin{equation}
        T_{n}(t) = \Phi_{n}(-\tau) \cos_{\tau}(\omega_{n}, t) + \tfrac{1}{\omega_{n}} \dot{\Phi}_{n}(-\tau) \sin_{\tau}(\omega_{n}, t) +
        \tfrac{1}{\omega_{n}} \int_{-\tau}^{0} \sin_{\tau}(\omega_{n}, t - \tau - s) \ddot{\Phi}(s) \mathrm{d}s.
        \label{EQ2_8}
    \end{equation}

    Thus, assuming sufficient smoothness of the data to be specified later,
    the solution $\xi_{0}$ to the homogeneous equation (\ref{EQ1_13}) satisfying homogeneous boundary and
    non-homogeneous initial conditions $\xi \equiv \Phi$ in $[-\tau, 0] \times [0, l]$ reads as
    \begin{equation}
        \begin{split}
            \xi_{0}(t, x) &= \sum_{n = 1}^{\infty} \left(\Phi_{n}(-\tau) \cos_{\tau}(\omega_{n}, t) +
            \tfrac{1}{\omega_{n}} \dot{\Phi}_{n}(-\tau) \sin_{\tau}(\omega_{n}, t) + \right. \\
            &\hspace{3cm} \left. \tfrac{1}{\omega_{n}} \int_{-\tau}^{0} \sin_{\tau}(\omega_{n}, t - \tau - s) \ddot{\Phi}_{n}(s) \mathrm{d}s\right)
            \sin\left(\tfrac{\pi n}{l} x\right), \\
            \Phi_{n}(t) &= \frac{2}{l} \int_{0}^{l} \left(\varphi(t, s) - G(t, s)\right) \sin\left(\tfrac{\pi n}{l} s\right) \mathrm{d}s
            \text{ for } n \in \NN.
        \end{split}
        \label{EQ2_9}
    \end{equation}

    \section{Non-homogeneous equation}
    Next, we consider the non-homogeneous equation (\ref{EQ1_16}) with
    the right-hand side from Equation (\ref{EQ1_18}) subject to homogeneous initial and boundary conditions
    \begin{equation}
        \partial_{tt} \xi_{1}(t, x) = a^{2} \partial_{xx} \xi_{1}(t - \tau, x) + c \xi_{1}(t - \tau, x) + F(t, x)
        \text{ for } (t, x) \in (0, T) \times (0, l), \notag
    \end{equation}
    where
    \begin{equation}
        F(t, x) = f(t, x) + c\left(\mu_{1}(t - \tau) + \tfrac{x}{l} \left(\mu_{2}(t - \tau) - \mu_{1}(t - \tau)\right)\right) -
        \left(\ddot{\mu}_{1}(t) - \tfrac{x}{l} \left(\ddot{\mu}_{2}(t) - \ddot{\mu}_{1}(t)\right)\right).
        \notag
    \end{equation}

    The solution will be constructed as as Fourier series with respect to the eigenfunctions of the Sturm \& Liouville problem
    from the previous section, i.e.,
    \begin{equation}
        \xi_{1}(t, x) = \sum_{n = 1}^{\infty} T_{n}(t) \sin\left(\tfrac{\pi n}{l} x\right).
        \label{EQ3_1}
    \end{equation}
    Plugging the ansatz from (\ref{EQ3_1}) into Equation (\ref{EQ1_6}) and comparing the time-dependent Fourier coefficients,
    we obtain a system of countably many second order delay differential equations
    \begin{equation}
        \ddot{T}_{n}(t) + \omega_{n} T_{n}(t - \tau) = F_{n}(t) \text{ for } t \in (0, T) \text{ with }
        F_{n}(t) = \frac{2}{l} \int_{0}^{l} F(t, s) \sin\left(\tfrac{\pi n}{l} x\right) \mathrm{d}s.
        \label{EQ3_2}
    \end{equation}

    In \cite{KhuDiRuLu2008}, the initial value problem for the non-homogeneous delay differential equation
    \begin{equation}
        \ddot{x}(t) + \omega^{2} x(t - \tau) = f(t) \text{ for } t \geq 0 \notag
    \end{equation}
    with homogeneous initial conditions $x \equiv 0$ in $[-\tau, 0]$ was shown to be uniquely solved by
    \begin{equation}
        x(t) = \int_{0}^{t} \sin_{\tau}(\omega, t - \tau - s) f(s) \mathrm{d}s.
        \label{EQ3_3}
    \end{equation}

    Exploiting Equation (\ref{EQ3_3}), the equations in (\ref{EQ3_2}) subject to zero initial conditions are uniquely solved by
    \begin{equation}
        T_{n}(t) = \int_{0}^{t} \sin_{\tau}(\omega_{n}, t - \tau - s) F_{n}(s) \mathrm{d}s.
        \label{EQ3_4}
    \end{equation}
    Therefore, the non-homogeneous partial delay differential equation with homogeneous initial and boundary condition formally reads as
    \begin{equation}
        \xi_{1}(t, x) = \sum_{n = 1}^{\infty} \left(\int_{0}^{t} \sin_{\tau}(\omega_{n}, t - \tau - s) F_{n}(s) \mathrm{d}s\right)
        \sin\left(\tfrac{\pi n}{l} x\right) \text{ for } (t, x) \in [0, T] \times [0, l].
        \label{EQ3_5}
    \end{equation}

    \section{General case solution}
    The solution in the general case can thus formally be represented as the following series
    \begin{equation}
        \begin{split}
            \xi(t, x) &= \sum_{n = 1}^{\infty} \left(\Phi_{n}(-\tau) \cos_{\tau}(\omega_{n}, t) + \tfrac{1}{\omega}_{n} \dot{\Phi}_{n}(-\tau) \sin_{\tau}(\omega_{n}, t) + \right. \\
            &\phantom{=\; \sum_{n = 1}^{\infty}} \left. \tfrac{1}{\omega_{n}} \int_{-\tau}^{0} \sin_{\tau}(\omega_{n}, t - \tau - s) \ddot{\Phi}_{n}(s) \mathrm{d}s\right)
            \sin\left(\tfrac{\pi n}{l} x\right) + \\
            &\phantom{=\;} \sum_{n = 1}^{\infty} \left(\int_{0}^{t} \sin_{\tau}(\omega_{n}, t - \tau - s) F_{n}(s) \mathrm{d}s\right)
            \sin\left(\tfrac{\pi n}{l} x\right) + G(t, x).
            \label{EQ3_6}
        \end{split}
    \end{equation}

    \subsection{Convergence of the Fourier series}

    \begin{theorem}
        Let $T > 0$, $\tau > 0$ and $m := \lceil\tfrac{T}{\tau}\rceil$.
        Further, let the data functions $\varphi$, $\mu_{1}$, $\mu_{2}$ and $f$
        be such that their Fourier coefficients $\Phi_{n}$ and $F_{n}$ given in Equations (\ref{EQ2_3}) and (\ref{EQ3_5})
        satisfy the conditions
        \begin{equation}
            \begin{split}
                \lim_{n \to \infty} \left(|\Phi_{n}(-\tau)| + |\dot{\Phi}_{n}(-\tau)|\right) n^{2m + 3 + \alpha} = 0, \quad
                \lim_{n \to \infty} \max_{s \in [-\tau, 0]} |\ddot{\Phi}_{n}| n^{2m + 3 + \alpha} = 0, \\
                \lim_{n \to \infty} \max_{k = 1, \dots, m} \max_{t \in [(k - 1) \tau, \max\{k \tau, T\}]}
                |F_{n}(t)| n^{2k + 3 + \alpha} = 0
            \end{split}
            \label{EQ3_7}
        \end{equation}
        for an arbitrary, but fixed $\alpha > 0$.
        Let
        \begin{equation}
            \left(\tfrac{\pi}{l} a\right)^{2} > c. \notag
        \end{equation}

        Then the classical solution to problem (\ref{EQ1_9})--(\ref{EQ1_11}) can be represented as an absolutely and uniformly convergent
        Fourier series given in Equation (\ref{EQ3_6}).
        The latter series is a twice continuously differentiable function with respect to both variables.
        Its derivatives of order less or equal two with respect to $t$ and $x$ can be obtained by a term-wise differentiation
        of the series and the resulting series are also absolutely and uniformly convergent in $[0, T] \times [0, l]$.
    \end{theorem}

    \begin{proof}
        We regroup the series from Equation (\ref{EQ3_6}) into the following sum
        \begin{equation}
            \xi(t, x) = S_{1}(t, x) + S_{2}(t, x) + S_{3}(t, x) + G(t, x), \notag
        \end{equation}
        where
        \begin{equation}
            \begin{split}
                S_{1}(t, x) &= \sum_{n = 1}^{\infty} A_{n}(t) \sin\left(\tfrac{\pi n}{l} x\right), \quad
                S_{2}(t, x) = \sum_{n = 1}^{\infty} B_{n}(t) \sin\left(\tfrac{\pi n}{l} x\right), \quad
                S_{3}(t, x) = \sum_{n = 1}^{\infty} C_{n}(t) \sin\left(\tfrac{\pi n}{l} x\right), \\
                A_{n}(t) &= \Phi_{n}(-\tau) \cos_{\tau}(\omega_{n}, t) + \tfrac{1}{\omega_{n}} \dot{\Phi}_{n}(-\tau) \sin_{\tau}(\omega_{n}, t), \\
                B_{n}(t) &= \tfrac{1}{\omega_{n}} \int_{-\tau}^{0} \sin_{\tau}(\omega_{n}, t - \tau - s) \ddot{\Phi}_{n} \mathrm{d}s, \\
                C_{n}(t) &= \tfrac{1}{\omega_{n}} \int_{-\tau}^{0} \sin_{\tau}(\omega_{n}, t - \tau - s) F_{n}(s) \mathrm{d}s
            \end{split}
            \notag
        \end{equation}
        and
        \begin{equation}
            \omega_{n} = \sqrt{\left(\tfrac{\pi n}{l} a\right)^{2} - c} \text{ for } n \in \NN. \notag
        \end{equation}
        \begin{enumerate}
            \item First, we consider the coefficient functions $A_{n}$. For an arbitrary $t \in [0, T]$ with
            $(k - 1) \tau \leq t < k\tau$, we find
            \begin{equation}
                \begin{split}
                    A_{n}(t) &= \Phi_{n}(-\tau) \cos_{\tau}(\omega_{n}, t) + \tfrac{1}{\omega_{n}} \dot{\Phi}_{n}(-\tau) \sin_{\tau}(\omega_{n}, t) \\
                    &= \left(1 - \left(\tfrac{\pi n}{l}a\right)^{2} \tfrac{t^{3}}{2!} + \dots +
                    (-1)^{k} \left(\tfrac{\pi n}{l}a\right)^{2k} \tfrac{(t - (k - 1) \tau)^{2k}}{(2k)!}\right) \Phi_{n}(-\tau) + \\
                    &\phantom{=;} \left((1 + \tau) - \left(\tfrac{\pi n}{l}a\right)^{2} \tfrac{t^{2}}{3!} + \dots +
                    (-1)^{k} \left(\tfrac{\pi n}{l}a\right)^{2k} \tfrac{(t - (k - 1) \tau)^{2k + 1}}{(2k + 1)!}\right) \Phi_{n}(-\tau).
                \end{split}
                \notag
            \end{equation}
            If $\Phi_{n}(-\tau)$ and $\dot{\Phi}_{n}(-\tau)$, $n \in \NN$, are such that the condition
            \begin{equation}
                \lim_{n \to \infty} \left(|\Phi_{n}(-\tau)| + |\dot{\Phi}_{n}(-\tau)|\right) n^{2k + 3 + \alpha} = 0 \notag
            \end{equation}
            holds true, the series $S_{1}$ as well as its derivatives of order less or equal 2 converge absolutely and uniformly.
            Note that a single differentiation with respect to $x$ corresponds, roughly speaking, to a multiplication with $n$.

            \item Next, we consider the coefficients $B_{n}$.
            For an arbitrary $t \in [0, T]$ with
            $(k - 1) \tau \leq t < k\tau$, we perform the substitution $t - \tau - s = \xi$ and exploit the mean value theorem to estimate
            \begin{equation}
                \begin{split}
                    |B_{n}(t)| &= \left|\tfrac{1}{\omega_{n}} \int_{t - \tau}^{t} \sin_{\tau}(\omega_{n}, \xi) \ddot{\Phi}_{n}(t - \tau - \xi) \mathrm{d}\xi\right| \\
                    &\leq \tau \max_{-\tau \leq s \leq 0} |\ddot{\Phi}_{n}(s)| \max_{j = k-1, k} \max_{t - \tau \leq s \leq t} \left|
                    (s - \tau) - \left(\tfrac{\pi n}{l}a\right)^2 \tfrac{s^{3}}{3!} + \dots \right. \\
                    &\hspace{5cm} \left. + (-1)^{j} \left(\tfrac{\pi n}{l}a\right)^{2j}
                    \tfrac{(s - (j - 1) \tau)^{2j+1}}{(2j + 1)!}\right|.
                \end{split}
                \notag
            \end{equation}
            Applying the theorem on differentiation under the integral sign to $B_{n}$ and taking into account that
            $\sin_{\tau}\left(\tfrac{\pi n}{l} a, \cdot\right)$ is twice weakly differentiable in $[0, \infty)$, namely: $\sin_{\tau}\left(\tfrac{\pi n}{l} a, \cdot\right) \in W^{2, \infty}_{\mathrm{loc}}\big((0, \infty)\big)$,
            its derivatives are polynomials of order lower than those of $\sin_{\tau}\left(\tfrac{\pi n}{l} a, \cdot\right)$
            and their convolution with $\ddot{\Phi}_{n}$ is continuous,
            analogous estimates can be obtained for $\dot{\Phi}_{n}$ and $\ddot{\Phi}_{n}$ which, in their turn,
            also follow to be continuous functions.

            Now, if the condition
            \begin{equation}
                \lim_{n \to \infty} \max_{s \in [-\tau, 0]} |\ddot{\Phi}_{n}| n^{2m + 3 + \alpha} = 0 \notag
            \end{equation}
            is satisfied, the series $S_{2}$ as well as its derivatives of order less or equal 2 converge absolutely and uniformly.

            \item Finally, we look at the Fourier coefficients $C_{n}$.
            Again, for an arbitrary period of time $t \in [0, T]$ with $(k - 1) t \leq t < k\tau$, $0 \leq k \leq m$,
            we substitute $t - \tau - \xi = s$.
            Once again, using the mean value theorem, we estimate
            \begin{equation}
                \begin{split}
                    |C_{n}(t)| &= \left|\tfrac{1}{\omega_{n}} \int_{t - \tau}^{t} \sin_{\tau}\left(\tfrac{\pi n}{l} a, \xi\right) F_{n}(t - \tau - \xi) \mathrm{d}\xi\right| \\
                    &\leq \tau \max_{t - \tau \leq s \leq t} |\ddot{\Phi}_{n}(s)| \max_{j = k - 1, k} \max_{t - \tau \leq s \leq t} \left|
                    (s - \tau) - \left(\tfrac{\pi n}{l} a\right)^{2} \tfrac{s^{3}}{3!} + \dots \right. \\
                    &\hspace{5cm} \left. + (-1)^{j} \left(\tfrac{\pi n}{l} a\right)^{2j} \tfrac{(s - (j - 1) \tau)^{2j+1}}{(2j + 1)!}\right|.
                \end{split}
                \notag
            \end{equation}
            As before, $C_{n}$ can be shown to be twice continuously differentiable.
            If now
            \begin{equation}
                \lim_{n \to \infty} \max_{k = 1, \dots, m} \max_{t \in [(k - 1) \tau, \max\{k \tau, T\}]}
                |F_{n}(t)| n^{2k + 3 + \alpha} = 0
            \end{equation}
            is satisfied, then both $S_{3}$ and its derivatives of order less or equal 2 converge absolutely and uniformly.
        \end{enumerate}
        Since all three conditions are guaranteed by the assumptions of the Theorem due to the fact $k \leq m$,
        the proof is finished.
    \end{proof}

    \begin{remark}
        From the practical point of view, the rapid decay condition on the Fourier coefficients of the data given in Equation (\ref{EQ3_7})
        mean a sufficiently high Sobolev regularity of the data and corresponding higher order
        compatibility conditions at the boundary of $(0, l)$ (cf. \cite{KhuPoAz2013}).
    \end{remark}

\end{document}